\documentclass[10pt]{article}
\usepackage{amsfonts,amssymb,amscd,amsmath,enumerate,verbatim,calc}

\title{\textbf{\Large{On some definite integrals connecting with certain infinite series}}}
\author{Ramesh Kumar Muthumalai
\footnote{Department of Mathematics, Sindhi college of arts and science, Chennai-77, Tamil Nadu, India.  Email id: ramjan\_80@yahoo.com.  Home page URL: http://ramjan\_07.page.tl/.}}
\date{}
\begin{document}
\maketitle
\paragraph{Abstract.} We show some definite integrals connecting to infinite series,  studied in Ramanujan's paper, titled "On question 330 of Professor Sanjana". We present few recursive methods to evaluate these definite integrals in various cases and we generalize this, to evaluate simliar kind of integrals through infinite series.
\paragraph{AMS subject classification:} 33E20; 33E50; 40A05.
\paragraph{Keywords:} Definite integrals; Question 330 of Professor Sanjana; Ramanujan papers; Infinite series.
\paragraph{1. Introduction.} The classical table of integrals by I.S. Gradshteyn and I.M. Ryzhik \cite{2} contains many entries related to definite integrals in the combinations of powers and algebraic functions of exponentials. In particular, the following integral is given for some special cases of $a$ and $\alpha$
\begin{equation}
\int_{0}^{\infty}x^{\alpha}e^{-bx}(1-e^{-x})^adx
\tag{1.1}
\end{equation}
In this paper, we connect this definite integral to infinite series given in \cite{3}. We evalute them for various cases $a$ and $\alpha$ by using recursive methods. Further, we have shown evaluation of combination of trigonometric functions and powers. Finally, we present summation of the series of this type 
\begin{equation}
\frac{1}{b^n}-\frac{a}{1!(b+1)^n}\frac{1}{a}+\frac{a(a-1)}{2!(b+2)^n}
\left(\frac{1}{a}+\frac{1}{a-1}\right)-\hdots
\tag{1.2}
\end{equation}
\paragraph{2. Definite integrals connecting with infinite series.} Consider the following definite integrals for $b>0$ and $\alpha\geq 0$
\begin{multline}
\int_0^\infty x^\alpha e^{-bx}(1-e^{-x})^adx \\
\qquad\quad\quad=\int_0^\infty x^\alpha e^{-bx} \left(1-\frac{a}{1!}e^{-x}+\frac{a(a-1)}{2!}e^{-2x}-\hdots\right)dx
\nonumber\\\qquad\qquad\quad=\Gamma(\alpha+1)\left(\frac{1}{b^{\alpha+1}}-
\frac{a}{1!(b+1)^{\alpha+1}}+\frac{a(a-1)}{2!(b+2)^{\alpha+1}}
-\hdots\right)\nonumber\\=\Gamma(\alpha+1)\phi(a,b,\alpha)
\quad\quad\qquad\qquad\qquad\qquad\qquad\qquad\qquad\quad\nonumber
\end{multline}
Where
\begin{equation}
\phi(a,b,\alpha)=\frac{1}{b^{\alpha+1}}-\frac{a}{1!(b+1)^{\alpha+1}}+
\frac{a(a-1)}{2!(b+2)^{\alpha+1}}-\hdots
\tag{2.1}
\end{equation}
Thus, we obtain
\begin{equation}
\int_0^\infty x^\alpha e^{-bx}(1-e^{-x})^adx =\Gamma(\alpha+1)\phi(a,b,\alpha)
\tag{2.2}
\end{equation}
If $\alpha=n$, $n\in N$ then
\begin{equation}
\int_0^\infty x^n e^{-bx}(1-e^{-x})^adx =n!\phi(a,b,n)
\tag{2.3}
\end{equation}
Also, (2.3) can be written as
\begin{equation}
\int_0^1 \log^n(1-t)(1-t)^{b-1}t^adt =(-1)^nn!\phi(a,b,n)
\tag{2.4}
\end{equation}
Differentiating (2.3) $m$ times with respect to $a$, we find that 
\begin{equation}
\int_0^\infty x^n e^{-bx}(1-e^{-x})^a\log^m(1-e^{-x})dx =n!\phi_a^{(m)}(a,b,n)
\tag{2.5}
\end{equation}
The above integral can be written as 
\begin{equation}
\int_0^1 \log^n(1-t)\log^m t(1-t)^{b-1}t^adt =(-1)^nn!\phi_a^{(m)}(a,b,n)
\tag{2.6}
\end{equation}
Also, it can be rewritten as 
\begin{equation*}
\int_0^1 \log^nt\log^m (1-t)t^{b-1}(1-t)^adt =(-1)^nn!\phi_a^{(m)}(a,b,n)
\end{equation*}
Now, using (2.6) and above equation we find
\begin{equation}
(-1)^mm!\phi_b^{(n)}(b-1,a+1,m) =(-1)^nn!\phi_a^{(m)}(a,b,n)
\tag{2.7}
\end{equation}
If $n=0$
\begin{equation}
(-1)^mm!\phi(b-1,a+1,m) =\phi_a^{(m)}(a,b,0)
\tag{2.8}
\end{equation}
Thus, through the infinite series $\phi(a,b,0)$, we can evaluate four different form of integrals. Similiarly, if we start from the infinite series for $b>0$ and $\alpha \geq 0$
\begin{equation}
\frac{1}{b^{\alpha+1}}+\frac{a}{1!(b+1)^{\alpha+1}}+
\frac{a(a-1)}{2!(b+2)^{\alpha+1}}+\hdots=\tilde{\phi}(a,b,\alpha)
\tag{2.9}
\end{equation}
then, we obtain
\begin{equation}
\int_0^\infty x^\alpha e^{-bx}(1+e^{-x})^adx =\Gamma(\alpha+1)\tilde{\phi}(a,b,\alpha)
\tag{2.10}
\end{equation}
Generalizing equations (2.2) and (2.10), for some  $|\beta|\leq 1$
\begin{equation}
\int_0^\infty x^\alpha e^{-bx}(1+\beta e^{-x})^adx =\Gamma(\alpha+1)\Psi(a,b,\beta,\alpha)
\tag{2.11}
\end{equation}
Where 
\begin{equation}
\frac{1}{b^{\alpha+1}}+\frac{a}{1!(b+1)^{\alpha+1}}\beta+
\frac{a(a-1)}{2!(b+2)^{\alpha+1}}\beta^2+\hdots=\Psi(a,b,\beta,\alpha)
\tag{2.12}
\end{equation}
Again, consider the following definite integral for  $b>0$ and $\alpha\geq 0$ 
\begin{multline}
\int_0^\infty x^\alpha e^{-ibx}(1-e^{-ix})^adx \\
\qquad\quad=\int_0^\infty x^\alpha e^{-ibx} \left(1-\frac{a}{1!}e^{-ix}+\frac{a(a-1)}{2!}e^{-2ix}-\hdots\right)dx
\nonumber\\
\qquad\qquad=\frac{\Gamma(\alpha+1)}{i^{\alpha+1}}\left(\frac{1}{b^{\alpha+1}}
-\frac{a}{1!(b+1)^{\alpha+1}}+\frac{a(a-1)}{2!(b+2)^{\alpha+1}}-\hdots
\right)\nonumber\\=-\Gamma(\alpha+1)\phi(a,b,\alpha)\left(\sin\frac{\alpha \pi}{2}+i\cos \frac{\alpha \pi}{2}\right)\quad\qquad\qquad\tag{2.13a}
\end{multline}
On the other hand, we have
\begin{align}
\int_0^\infty x^\alpha e^{-ibx}(1-e^{-ix})^adx= \qquad\qquad\qquad\qquad\qquad\qquad\qquad\qquad\nonumber\\
 2^a\int_0^\infty x^\alpha \sin^\alpha\frac{x}{2}\left(\cos \left[\frac{\pi-x}{2}a+bx\right]-i\sin \left[\frac{\pi-x}{2}a+bx\right]\right)dx
\tag{2.13b}
\end{align}
Comparing (2.13a) and (2.13b), we obtain 
\begin{align}
\int_0^\infty x^\alpha \sin^a \frac{x}{2}\cos \left(\frac{\pi-x}{2}a+bx\right)dx=-2^{-a}\Gamma(\alpha+1)\phi(a,b,\alpha)
\sin\frac{\alpha \pi}{2} 
\tag{2.14}
\end{align}
\begin{align}
\int_0^\infty x^\alpha \sin^a \frac{x}{2}\sin \left(\frac{\pi-x}{2}a+bx\right)dx=2^{-a}\Gamma(\alpha+1)\phi(a,b,\alpha)
\cos\frac{\alpha \pi}{2} 
\tag{2.15}
\end{align}
Expanding the  cosine funtion in the integral (2.14),
\begin{align}
\int_0^\infty x^\alpha \sin^a \frac{x}{2}\left(\cos \frac{\pi a}{2}\cos\left(b-\frac{a}{2}\right)x-\sin \frac{\pi a}{2}\sin\left(b-\frac{a}{2}\right)x\right)\nonumber\\
=-2^{-a-\alpha-1}\Gamma(\alpha+1)\phi(a,b,\alpha)\sin\frac{\alpha \pi}{2} 
\nonumber
\end{align}
Setting $u=2b-a$, then (2.14) gives
\begin{align}
\int_0^\infty x^\alpha \sin^a x\left(\cos \frac{\pi a}{2}\cos ux-\sin \frac{\pi a}{2}\sin ux\right)\qquad\qquad\nonumber\\
=-2^{-a-\alpha-1}\Gamma(\alpha+1)\phi\left(a,\frac{u+a}{2},\alpha\right)\sin
\frac{\alpha \pi}{2} 
\nonumber
\end{align}
Setting $\lambda_c=\int_0^\infty x^\alpha \sin^a x \cos uxdx$ and $\lambda_s=\int_0^\infty x^\alpha \sin^a x \sin uxdx$
\begin{align}
\cos \frac{\pi a}{2}\lambda_c-\sin \frac{\pi a}{2}\lambda_s
=-2^{-a-\alpha-1}\Gamma(\alpha+1)\phi\left(a,\frac{u+a}{2},\alpha\right)\sin
\frac{\alpha \pi}{2} 
\tag{2.16a}
\end{align}
Similialry from (2.15), we have
\begin{align}
\sin \frac{\pi a}{2}\lambda_c+\cos \frac{\pi a}{2}\lambda_s
=2^{-a-\alpha-1}\Gamma(\alpha+1)\phi\left(a,\frac{u+a}{2},\alpha\right)\cos
\frac{\alpha \pi}{2} 
\tag{2.16b}
\end{align}
Solving (2.16a) and (2.16b), we get
\begin{equation}
\lambda_c=2^{-a-\alpha-1}\Gamma(\alpha+1)\phi\left(a,\frac{u+a}{2},\alpha\right)
\sin(a-\alpha)\frac{ \pi}{2} 
\tag{2.17}
\end{equation}
\begin{equation}
\lambda_s
=2^{-a-\alpha-1}\Gamma(\alpha+1)\phi\left(a,\frac{u+a}{2},\alpha\right)
\cos(a-\alpha)\frac{ \pi}{2} 
\tag{2.18}
\end{equation}
Differentiating (2.17) and (2.18) $m$ times, with respect to $a$, we can evaluate the following definite integrals
\begin{equation}
\int_0^\infty x^\alpha  \log^m\sin x\sin^a x \cos uxdx \hspace{0.3cm}\mbox{and} \hspace{0.3cm}\int_0^\infty x^\alpha  \log^m\sin x\sin^a x \sin uxdx
\nonumber
\end{equation}
Similiarly, if we start from the integral $\int_0^\infty x^\alpha e^{-ibx}(1+e^{-ix})^adx$ we find
\begin{align}
\int_0^\infty x^\alpha \cos^a \frac{x}{2}\cos \left(\frac{a}{2}+b\right)xdx=-2^{-a}\Gamma(\alpha+1)\tilde{\phi}(a,b,\alpha)
\sin\frac{\alpha \pi}{2} 
\nonumber
\end{align}
\begin{align}
\int_0^\infty x^\alpha \cos^a \frac{x}{2}\sin \left(\frac{a}{2}+b\right)xdx=2^{-a}\Gamma(\alpha+1)\tilde\phi(a,b,\alpha)
\cos\frac{\alpha \pi}{2} 
\nonumber
\end{align}
Setting $v=a+2b$ and after simplification, we find
\begin{equation}
\int_0^\infty x^\alpha \cos^a x\cos vxdx =-2^{-a-\alpha-1}\Gamma(\alpha+1)\tilde{\phi}\left(a,\frac{v-a}{2},\alpha\right)
\sin\frac{\alpha \pi}{2} 
\tag{2.19}
\end{equation}
\begin{equation}
\int_0^\infty x^\alpha \cos^a x\sin vxdx=2^{-a-\alpha-1}\Gamma(\alpha+1)\tilde\phi \left(a,\frac{v-a}{2},\alpha\right)\cos\frac{\alpha \pi}{2} 
\tag{2.20}
\end{equation}
\paragraph{3. Algebraic recursive method.} Let us start by deriving a recursive method connecting algebraic functions to infinite series, for  some $0<|t|<1$, we know that from binomial theorem [1] 
\begin{align}
1-\frac{p}{1!}t+\frac{p(p-1)}{2!}t^2-\hdots=&(1-t)^p\nonumber\\
=&A_1^{(1)}(1-t)^p
\nonumber
\end{align}
Where $A_1^{(1)}=1$ and 
\begin{align}
b-\frac{p}{1!}(b+1)t+\frac{p(p-1)}{2!}(b+2)t^2-\hdots=&b(1-t)^p-pt(1-t)^{p-1}
\nonumber\\=&A_1^{(2)}(1-t)^p+A_2^{(2)}t(1-t)^{p-1}\nonumber
\end{align}
Where $A_1^{(2)}=b$ and $A_2^{(2)}=-p$. Similiarly, we can write
\begin{align}
b^2-\frac{p}{1!}(b+1)^2t+\frac{p(p-1)}{2!}(b+2)^2t^2-\hdots
\qquad\qquad\qquad\qquad
\nonumber\\\qquad\qquad\qquad
=b^2(1-t)^p-(2b+1)pt(1-t)^{p-1}+p(p-1)t^2(1-t)^{p-2}
\nonumber \\= A_1^{(3)}(1-t)^p+A_2^{(3)}t(1-t)^{p-1}+A_3^{(3)}t^2(1-t)^{p-2}
\nonumber\qquad\quad
\end{align}
Where $A_1^{(3)}=b^2$, $A_2^{(3)}=-p(2b+1)$ and $A_3^{(3)}=p(p-1)$.
Generalizing  this to, some positive integer $m$, we have
\begin{align}
\sum_{i=0}^\infty(-1)^i\binom{p}{i}(b+i)^{m-1}t^i=
\sum_{k=1}^mA_k^{(m)}t^{k-1}(1-t)^{p-k+1}
\tag{3.1}
\end{align}
Where $A^{(m)}_k\hspace{0.2cm}(k=1,2,,\hdots,m$ and $m=1,2,3,\hdots)$ are  independent of $t$. To determine all  other  unknown $A's$, multiply both sides of (3.1) by $t^b$ 
\begin{align}
\sum_{i=0}^\infty(-1)^i\binom{p}{i}(b+i)^{m-1}t^{i+b}&=
\sum_{k=1}^mA_k^{(m)}t^{k+b-1}(1-t)^{p-k+1}
\nonumber
\end{align}
Differentiating with respect to $t$,
\begin{align}
\sum_{i=0}^\infty(-1)^i\binom{p}{i}(b+i)^{m}t^{i+b-1}\qquad\qquad\qquad\qquad
\qquad\qquad\qquad\qquad\qquad\qquad\qquad\nonumber\\
=\sum_{k=1}^mA_k^{(m)}
\left((k+b-1)t^{k+b-2}(1-t)^{p-k+1}-(p-k+1)t^{k+b-1}(1-t)^{p-k}\right)
\nonumber
\end{align}
Now rearranging
\begin{align}
=A_1^{(m)}bt^{b-1}(1-t)^p\qquad\qquad\qquad\qquad\qquad\qquad
\qquad\qquad\qquad\qquad\qquad\nonumber\\
+\sum_{k=2}^m
\left(-(p-k+2)A_{k-1}^{(m)}+(b+k-1)A_k^{(m)}\right)t^{b+k-2}(1-t)^{p-k+1}
\nonumber\\-(p-m+1)A_m^{(m)}t^{m+b-1}(1-t)^{p-m}\nonumber
\end{align}
Cancelling $t^{b-1}$ on both sides
\begin{align}
A_1^{(m)}b(1-t)^p\qquad\qquad\qquad\qquad\qquad\qquad
\qquad\qquad\qquad\qquad\qquad\nonumber\\+\sum_{k=2}^m\left(-(p-k+2)
A_{k-1}^{(m)}+(b+k-1)A_k^{(m)}\right)t^{k-1}(1-t)^{p-k+1}\nonumber\\
-(p-m+1)A_m^{(m)}t^{m}(1-t)^{p-m}=\sum_{i=0}^\infty(-1)^i\binom{p}{i}
(b+i)^{m}t^i
\tag{3.2}
\end{align}
Now (3.2) can be written in the form of equation (3.1). If we replace $m$ by $m+1$ in (3.1) 
\begin{align}
\sum_{k=1}^{m+1}A_k^{(m+1)}t^{k-1}(1-t)^{p-k+1}=\sum_{i=0}^\infty(-1)^i
\binom{p}{i}(b+i)^{m}t^i
\tag{3.3}
\end{align}
Comparing (3.2) and (3.3)
\begin{align} 
A_1^{(m+1)}=&bA_1^{(m)}\nonumber\\
A_k^{(m+1)}=&-(p-k+2)A_{k-1}^{(m)}+(b+k-1)A_k^{(m)}, \hspace{1cm} k=2,3,\hdots m \tag{3.4}\\
A_{m+1}^{(m+1)}=&-(p-m+1)A_m^{(m)}\nonumber
\end{align}
Hence, we can evaluate all A's recursively through equation (3.4).\\\\
\textbf{4. Evaluation of definite integrals through infinite series.} In section 2, we have seen some definite integrals connecting with infinite series. Now, we will study evalutaion of these infinite series for various cases.\\\\
\textbf{Evaluation of $\phi(a,b,\alpha)$}\\\\
The integral $\int_0^\infty x^\alpha e^{-bx}(1-e^{-x})^adx$ exists, if  $\alpha \geq 0$ and $b>0$. It is clear that, from equation (2.2)
\begin{equation}
\int_0^\infty x^\alpha e^{-bx}(1-e^{-x})^adx =\Gamma(\alpha+1)\phi(a,b,\alpha)
\nonumber
\end{equation}
Now, we evaluate them for various cases of $\alpha$ and $a$ through infinite series $\phi(a,b,\alpha)$.\\\\
\textbf{Case 1.} If $\alpha=n$ is an integer, then  we can use Ramanujan formula, to evaluate infinite series $\phi(a,b,\alpha)$ as follow as 
\begin{equation}
n\phi(a,b,,n)=\sigma_1\phi(a,b,,n-1)+\sigma_2\phi(a,b,,n-2)+\hdots+
\sigma_n\phi(a,b,,0)
\tag{4.1}
\end{equation}
Where
\begin{equation}
\sigma_k=\frac{1}{b^k}-\frac{1}{(a+b+1)^k}+\frac{1}{(b+1)^k}-\frac{1}{(a+b+2)^k}
+\hdots 
\tag{4.2}
\end{equation}
and $k=1,2,3,\hdots$
\\\\\textbf{Example.} If we take $a=-\frac{1}{2}$, $b=\frac{1}{4}$ and $n=0,1,2$, then $\phi(-\frac{1}{2},\frac{1}{4},n)$ can be calculated from (4.1). Finally, using (2.2) we find
\begin{equation*}
\int_0^\infty x e^{-\frac{x}{4}}(1-e^{-x})^{-\frac{1}{2}}dx =\frac{\Gamma(\frac{1}{4})^2\pi}{\sqrt{2\pi}}\qquad\qquad\qquad\qquad
\end{equation*}
\begin{equation*}
\int_0^\infty x^2 e^{-\frac{x}{4}}(1-e^{-x})^{-\frac{1}{2}}dx=\frac{\Gamma(\frac{1}{4})^2}
{\sqrt{2\pi}}\left(\pi^2+16S'_2\right)\qquad\qquad
\end{equation*}
\begin{equation*}
\int_0^\infty x^3 e^{-\frac{x}{4}}(1-e^{-x})^{-\frac{1}{2}}dx =\frac{\Gamma(\frac{1}{4})^2}{\sqrt{2\pi}}\left(5\pi^3+48S'_2+128S'_3\right)
\end{equation*}
Where
\begin{equation*}
S'_r=\frac{1}{1^r}-\frac{1}{3^r}+\frac{1}{5^r}-\frac{1}{7^r}+\hdots+\infty 
\end{equation*}\\
\textbf{Case 2.} If $\alpha=\mu+m, \mu>0, m\in N$ and $a=-1,-2,-3,\hdots$ then putting $p=-1, t=e^{-x}$ in (3.3) and multiplying by $x^{\mu+m}e^{-bx}$ on both sides, we find
\begin{equation*}
\sum_{k=1}^{m+1}A_{k}^{(m+1)}\frac{x^{\mu+m}e^{-x(b+k-1)}}{(1-e^{-x})^k}
=\sum_{i=0}^\infty(b+i)^mx^{\mu+m}e^{-(b+i)x}
\end{equation*}
Integrating on $(0,\infty)$,
\begin{equation}
\sum_{k=1}^{m+1}A_{k}^{(m+1)}\int_0^\infty\frac{x^{\mu+m}e^{-x(b+k-1)}}
{(1-e^{-x})^k}dx=\sum_{i=0}^\infty(b+i)^m\int_0^\infty x^{\mu+m}e^{-(b+i)x}dx
\tag{4.3}
\end{equation}
Using (2.2) and after simplification, we get
\begin{equation}
\sum_{k=1}^{m+1}A_{k}^{(m+1)}\phi(-k,b+k-1\mu+m)=\zeta(\mu+1,b)
\tag{4.4}
\end{equation}
\textbf{Example.} If we put $m=0$ and $\mu>0$ in (4.4)
\begin{equation*}
\phi(-1,b,\mu)=\zeta(\mu+1,b)
\end{equation*}
Putting $m=1$ in (4.4)
\begin{equation*}
\phi(-2,b+1,\mu+1)=\zeta(\mu+1,b)-b\zeta(\mu+2,b)
\end{equation*}
Putting $m=2$ in (4.4)
\begin{equation*}
\phi(-3,b+2,\mu+2)=\zeta(\mu+1,b)-(2b+1)\zeta(\mu+2,b)+b(b+1)\zeta(\mu+3,b)
\end{equation*}
\textbf{Case 3.} If $\alpha=\mu+m, \mu>0, m\in N$ and $a=\gamma+m\neq-1,-2,-3,\hdots$ then putting $\gamma+m=p, t=e^{-x}$ in (3.3) and multiplying by $x^{\mu+m}e^{-bx}$ on both sides and integrating on $(0,\infty)$
\begin{align*}
\sum_{k=1}^{m+1}A_{k}^{(m+1)}\int_0^\infty x^{\mu+m}e^{x(-b-k+1)}{(1-e^{-x})}^{\gamma+m-k+1}
dx\qquad\qquad\qquad\qquad\nonumber\\=\sum_{i=0}^\infty(-1)^i\binom{\gamma+m}
{i}(b+i)^m\int_0^\infty x^{\mu+m}e^{-(b+i)x}dx
\end{align*}
Using (2.2) and after simplification, we get
\begin{equation}
\sum_{k=1}^{m+1}A_{k}^{(m+1)}\phi(\gamma+m-k+1,b+k-1,\mu+m)=\phi(\gamma+m,b,\mu)
\tag{4.5}
\end{equation}
\textbf{Example.} If we take $\gamma=-\frac{1}{2}$, $b=\frac{1}{4}$, $\mu=0$ and $m=1$ in (3.3), then
\begin{equation*}
\frac{1}{4}\phi\left(-\frac{1}{2},\frac{1}{4},1\right)+\frac{1}{2}\phi
\left(-\frac{3}{2},\frac{5}{4},1\right)
=\phi\left(-\frac{1}{2},\frac{1}{4},0\right)
\end{equation*}
So that
\begin{equation*}
\int_0^\infty x e^{-\frac{5x}{4}}(1-e^{-x})^{-\frac{3}{2}}dx =2\frac{\Gamma(\frac{1}{4})^2}{\sqrt{2\pi}}\left(1-\frac{\pi}{4}\right)
\qquad\qquad\qquad\qquad
\end{equation*}
\textbf{Evaluation of $\tilde{\phi}(a,b,\alpha)$}\\\\
This integral exists, if  $\alpha \geq 0$ and $b>0$. It is clear that, from the equation (2.10).  
\begin{equation}
\int_0^\infty x^\alpha e^{-bx}(1+e^{-x})^adx =\Gamma(\alpha+1)\tilde{\phi}(a,b,\alpha)
\nonumber
\end{equation}
Now, we evaluate them for various cases of  $\tilde{\phi}(a,b,\alpha)$.\\\\
\textbf{Case 1.}If $\alpha=\mu+m, \mu \geq 0, m\in N$ and $a=-1,-2,-3,\hdots$\\\\
If $a=1$, then
\begin{equation}
\int_0^\infty \frac{x^\alpha e^{-bx}}{1-e^{-2x}}dx =\frac{\Gamma(\alpha+1)}{2^\alpha}\phi\left(-1,\frac{b}{2},\alpha\right)
\nonumber
\end{equation}
Also, 
\begin{equation}
\int_0^\infty \frac{x^\alpha e^{-bx}}{1+e^{-x}}dx =\frac{\Gamma(\alpha+1)}{2^\alpha}\phi\left(-1,\frac{b}{2},\alpha\right)-
\Gamma(\alpha+1)\phi(-1,b,\alpha)
\nonumber
\end{equation}
\begin{equation}
\tilde{\phi}(-1,b,\alpha)=\frac{1}{2^\alpha}\phi
\left(-1,\frac{b}{2},\alpha\right)-\phi(-1,b,\alpha)
\tag{4.6}
\end{equation}
Putting $p=-1, t=e^{-x}$ in (3.3), multiplying by $x^{\mu+m}e^{-bx}$ on both sides and integrating $(0,\infty)$, we get 
\begin{equation}
\sum_{k=1}^{m+1}A_{k}^{(m+1)}(-1)^{k-1}\tilde{\phi}(-k,b+k-1,\mu+m)
=\sum_{i=0}^\infty\frac{(-1)^i}{(b+i)^{\mu+1}}
\tag{4.7}
\end{equation}
\textbf{Case 2.} If $\alpha=\mu+m, \mu>0$ and $a=\gamma+m\neq-1,-2,-3,\hdots$. Simliliarly, putting $p=\gamma+m, t=-e^{-x}$ in (3.5),  multiplying by $x^{\mu+m}e^{-bx}$ on both sides and integrating on $(0,\infty)$
\begin{equation}
\sum_{k=1}^{m+1}(-1)^{k-1}A_{k}^{(m+1)}\tilde{\phi}(\gamma+m-k+1,b+k-1,\mu+m)
=\tilde{\phi}(\gamma+m,b,\mu)
\tag{4.8}
\end{equation}\\
\textbf{Evaluation of $\Psi(a,b,\beta,\alpha)$}\\\\
Putting $p=\gamma+m,  t=-\beta e^{-x}$,  $|\beta|\leq 1$ in (3.3),  multiplying by $x^{\mu+m}e^{-bx}$ and integrating on $(0,\infty)$, then
\begin{equation}
\sum_{k=1}^{m+1}A_{k}^{(m+1)}(-\beta)^{k-1}\Psi(\gamma+m-k+1,b+k-1,\beta,\mu+m)=
\Psi(\gamma+m,b,\beta,\mu)
\tag{4.9}
\end{equation}
\textbf{Example.}
If we put $m=0$ and $\mu\geq 0$ in (4.9)
\begin{equation*}
\Psi(-1,b,\beta,\mu)=\Phi(\beta,\mu+1,b)
\end{equation*}
Putting $m=1$ in (4.9)
\begin{equation*}
\Psi(-2,b+1,\beta, \mu+1)=\frac{1}{\beta}\left(\Phi(\beta,\mu+1,b)-b\Phi(\beta,\mu+2,b)\right)
\end{equation*}
For another interesting example, consider the following definite integral from \cite {2}, pp 350.
\begin{equation}
\int_{-\infty}^\infty \frac{x^2 e^{-bx}}{1+e^{-x}}dx =\pi^3 \csc b\pi(2-\sin^2b\pi)
\nonumber
\end{equation}
Now, for $p=-1,-2,-3,\hdots$, $t=-\beta e^{-x}$ in (3.5), multiply and divide by $x^2e^{-bx}$ and $1+\beta e^{-x}$ respectively, then, after integrating on  $(-\infty,\infty)$
\begin{equation}
\sum_{k=1}^{m+1}A^{(m+1)}_k\int_{-\infty}^\infty \frac{(-\beta)^{k-1}x^2 e^{-x(b+k-1)}}{(1+e^{-x})(1+\beta e^{-x})^k}dx =\sum_{i=0}^{\infty}\beta^i(b+i)^m\int_{-\infty}^\infty \frac{x^2 e^{-(b+i)x}}{1+e^{-x}}dx
\nonumber
\end{equation}
\begin{equation}
 =\pi^3 \csc b\pi(2-\sin^2b\pi)\sum_{i=0}^{\infty}\beta^i(b+i)^m
\nonumber
\end{equation}
If we put $m=0$, we obtain
\begin{equation}
\int_{-\infty}^\infty \frac{x^2 e^{-bx}}{(1+e^{-x})(1+\beta e^{-x})}dx =
\frac{\pi^3}{1-\beta} \csc b\pi(2-\sin^2b\pi)
\nonumber
\end{equation}
By putting $m=1,2,\hdots$, we can easily evaluate definite  integrals of the type \\$\int_{-\infty}^\infty \frac{x^2 e^{-bx}}{(1+e^{-x})(1+\beta e^{-x})^{m+1}}dx $.
Thus, using new algebraic recursive method, we can evaluate simlilar type of integrals.\\\\
\textbf{5. Summation of Infinite series.} To evaluate the definite integrals  given in equation (2.6), we  make use of summation of the series of the type (1.2) by differentiating Ramanujan recursive formula (2.4) sucessively with respect to $a$. but, it may require some additional calculation 
 and is time consuming. Instead, we can use the below formula which need less computations then the differentiation of (2.4).  We derive this formula almost similiar way of Ramanujan formula. For $a \neq -1,-2,-3,\hdots$, using Taylor series [2] near $p=0$
\begin{equation}
\log\Gamma(p+a+1)=\log\Gamma(a+1)+\psi(a+1)\frac{p}{1}+\zeta(2,a+1)\frac{p^2}{2} -\zeta(3,a+1)\frac{p^3}{3}+\hdots
\tag{5.1}
\end{equation}
Where
\begin{equation}
\psi(a+1)=\lim_{n \to \infty}\log n-\left(\frac{1}{a+1}+\frac{1}{a+2}+\frac{1}{a+3}+\hdots+\frac{1}{a+n}\right)
\tag{5.2}
\end{equation}
and
\begin{equation}
\zeta(r,a+1)=\frac{1}{(a+1)^r}+\frac{1}{(a+2)^r}+\frac{1}{(a+3)^r}+\hdots
\tag{5.3}
\end{equation}
Let us take 
\begin{equation*}
f(p,a,b)=\int_0^1x^{p+b-1}(1-x)^adx
\end{equation*}
Expanding in ascending powers of $p$ and integrating by $x$,
\begin{equation*}
=\frac{1}{b+p}-\frac{a}{1!(b+p+1)}+\frac{a(a-1)}{2!(b+p+2)}
-\frac{a(a-1)(a-2)}{3!(b+p+3)}+\hdots
\end{equation*}
\begin{equation}
f(p,a,b)=\phi(a,b,0)-p\phi(a,b,1)+p^2\phi(a,b,2)-p^3\phi(a,b,3)+\hdots
\tag{5.4}
\end{equation}
Now, differentiating with respect to $a$,
\begin{equation}
f'(p,a,b)=\phi'_a(a,b,0)-p\phi'_a(a,b,1)+p^2\phi'_a(a,b,2)
+\hdots
\tag{5.4a}
\end{equation}
\begin{equation*}
f(p,a,b)=\frac{\Gamma(p+b)\Gamma(a+1)}{\Gamma(p+a+b+1)}
\end{equation*}
Therefore, 
\begin{equation*}
\log f(p,a,b)=\log \Gamma(p+b) +\log \Gamma(a+1)-\log \Gamma(p+a+b+1)
\end{equation*}
Again, differentiating with respect to $a$
\begin{equation*}
\frac{f'(p,a,b)}{f(p,a,b)}=\psi(a+1) -\psi(p+a+b+1)
\end{equation*}
\begin{equation*}
\frac{f'(p,a,b)}{f(p,a,b)}=\psi(a+1) -\left(\psi(a+b+1)+p\zeta(2,a+b+1)-p^2\zeta(3,a+b+1)+\hdots\right)
\end{equation*}
Using (5.4) and (5.4a)
\begin{align*}
\phi'_a(a,b,0)-p\phi'_a(a,b,1)+p^2\phi'_a(a,b,2)-\hdots
\qquad\qquad\qquad\qquad\nonumber\\=\left(\psi(a+1) -\psi(a+b+1)-p\zeta(2,a+b+1)+p^2\zeta(3,a+b+1)-\hdots\right)\times\\
\left(\phi(a,b,0)-p\phi(a,b,1)+p^2\phi(a,b,2)-p^3\phi(a,b,3)+\hdots\right)
\end{align*}
Equating coefficients of $p$
\begin{equation*}
\phi'_a(a,b,0)=(\psi(a+1) -\psi(a+b+1))\phi(a,b,0)
\end{equation*}
\begin{equation*}
\phi'_a(a,b,1)=(\psi(a+1) -\psi(a+b+1))\phi(a,b,1)+\zeta(2,a+b+1)\phi(a,b,0)
\end{equation*}
\begin{align*}
\phi'_a(a,b,2)=(\psi(a+1) -\psi(a+b+1))\phi(a,b,2)\nonumber\qquad\qquad\qquad\qquad\qquad\\
+\zeta(2,a+b+1)\phi(a,b,1)+\zeta(3,a+b+1)\phi(a,b,0)
\end{align*}
Generalizing this, we obtain
\begin{align}
\phi'_a(a,b,n)=(\psi(a+1)-\psi(a+b+1))\phi(a,b,n)\qquad\qquad\qquad\nonumber\\
+\zeta(2,a+b+1)\phi(a,b,n-1)+\zeta(3,a+b+1)\phi(a,b,n-2)\nonumber\\
+\hdots+\zeta(n+1,a+b+1)\phi(a,b,0)
\tag{5.5}
\end{align}
Further $\phi''_a(a,b,n),\phi'''_a(a,b,n),\hdots$ can be found by successive differentiation of (5.5). Also, summation of the series (1.2) is given by 
\begin{align}
\frac{1}{b^n}-\frac{a}{1!(b+1)^n}\frac{1}{a}+\frac{a(a-1)}{2!(b+2)^n}\left(
\frac{1}{a}+\frac{1}{a-1}\right)-\hdots\nonumber\\=\frac{1}{b^n}
+\phi'_a(a,b,n-1)
\tag{5.6}
\end{align}
Simplifying (5.6) and letting $a$ tends to 0, we obtain
\begin{equation}
\frac{1}{1(b+1)^n}+\frac{1}{2(b+2)^n}+\frac{1}{3(b+3)^n}-
\hdots=\phi'_a(0,b,n-1)
\tag{5.7}
\end{equation}
Thus, from (5.5)
\begin{equation}
\phi'_a(0,b,n)=-\frac{C+\psi(b+1)}{b^{n+1}}+\frac{\zeta(2,b+1)}{b^n}+
\hdots+\frac{\zeta(n+1,b+1)}{b}
\tag{5.8}
\end{equation}
Where $C$ is  Euler constant.\\\\
\textbf{Examples.} If we take $a=-\frac{1}{2}$ and $b=\frac{1}{4}$ then
\begin{align*}
\frac{1}{1^n}+\frac{1}{2.5^n}\frac{1}{1}+\frac{1.3}{2.4.9^n}\left(\frac{1}{1}+
\frac{1}{3}\right)+\frac{1.3.5}{2.4.6.13^n}\left(\frac{1}{1}+\frac{1}{3}
+\frac{1}{5}\right)\hdots\\
=1-\frac{1}{2.4^n}\phi'_a\left(-\frac{1}{2},\frac{1}{4},n-1\right)
\end{align*}
If $n=0,1$, then $\phi'\left(-\frac{1}{2},\frac{1}{4},n\right)$ can be calculated from (5.5). Finally, using (2.6) we find
\begin{equation*}
\int_0^1 \log t (1-t)^{-\frac{3}{4}}t^{-\frac{1}{2}}dt =\frac {\Gamma\left(\frac{1}{4}\right)^2}{\sqrt{2\pi}}\left(-\frac{\pi}{2}+\ln 2\right)
\end{equation*}
\begin{equation*}
\int_0^1 \log(1-t)\log t (1-t)^{-\frac{3}{4}}t^{-\frac{1}{2}}dt =\frac {\Gamma\left(\frac{1}{4}\right)^2}{\sqrt{2\pi}}\left(\pi\left(-\frac{\pi}{2}+\ln 2\right)+\zeta\left(2,\frac{3}{4}\right)\right)
\end{equation*}\\
\textbf{6. Conclusion.} In conclusion we note that the new recursive formulas to evaluate  definite integrals in the combination of powers with algebraic functions of exponentials and  trigonometric functions has been developed in this article. First, we have studied some  definite integrals connecting with  Ramanujan \cite{3} infinite series, Riemann zeta function, etc and then we derive various  recursive formulas to evaluate these infinite series. 

\end{document}